\documentclass[a4paper,conference]{IEEEtran}
\IEEEoverridecommandlockouts

\usepackage{cite,color,url}
\usepackage{amsmath,amssymb,amsfonts,hyperref,multirow}
\usepackage{booktabs}
\usepackage{algorithmic,tikz}
\usepackage{subcaption,amsthm,longtable,lscape,amsmath}
\usepackage{graphicx,longtable}
\usepackage{textcomp}
\usepackage{xcolor}
\usetikzlibrary{fit,calc}
\def\BibTeX{{\rm B\kern-.05em{\sc i\kern-.025em b}\kern-.08em
    T\kern-.1667em\lower.7ex\hbox{E}\kern-.125emX}}
\begin{document}

\title{On Solving the Set Covering Problem \\with Conflicts on Sets}
\author{\IEEEauthorblockN{Roberto Montemanni}
\IEEEauthorblockA{\textit{Department of Sciences and Methods for Engineering} \\
\textit{University of Modena and Reggio Emilia}\\
Reggio Emilia, Italy \\
roberto.montemanni@unimore.it}
\and
\IEEEauthorblockN{Derek H. Smith}
\IEEEauthorblockA{\textit{Faculty of Computing, Engineering and Science} \\
\textit{University of South Wales}\\
Pontypridd, Wales, UK \\
derek.smith@southwales.ac.uk}}

\maketitle
\begin{abstract}
A variant of the well-known Set Covering Problem is studied in this paper, where subsets of a collection have to be selected, and pairwise conflicts among subsets of items exist. The selection of each subset has a cost, and the inclusion of conflicting subsets is associated with a penalty to be paid. The problem, which can be used to model real applications, looks for a selection of subsets that  cover the original collection, while
minimizing the sum of covering and penalty costs. 

In this paper we consider a compact mixed integer linear
program  and we solve it with an open-source solver.
Computational results on the benchmark instances commonly
used in the literature of the problem are reported. The results
indicate that the new approach we propose is capable of good results, both in terms of lower and upper bounds, although not matching the state-of-the-art on average. The new approach was, however, able to improve 9 best-known heuristic solutions.
\end{abstract}


\section{Introduction}
The Set Covering Problem (SCP) is about fully covering  a set of elements with a collection of subsets of such elements, where each subset is characterized by a given cost, and the objective is to find the feasible solution with minimum cost.

A generalization of the SCP was introduced in Carrabs et al. \cite{car24}, motivated by real applications. Some pairs of subsets can be \emph{in conflict}, and a penalty is incurred if both the subsets of the pair are selected at the same time.The new problem is named the Set Covering Problem with Conflicts on Sets (SCP-CS). Similar problems can be found in \cite{jac19} named the \emph{Conflict-free Set Covering Problem},  where a variant is considered with pairs of subsets that are classified as incompatible and cannot therefore be selected simultaneously.  The study is focussed on approximation algorithms and computational complexity. Different families of valid inequalities for a model representing  the same problem were presented in \cite{saf22}, together with  some preprocessing routines introduced to reduce the sizes of the instances and procedures to further speed up the problem
resolution.  A specialization of the same problem in geometric settings was discussed in \cite{ban20}, where some approximation results were presented.

As suggested in \cite{car24}, an actual application related to the SCP-CS problem can be identified in telecommunications. Modern society is increasingly relying on wireless data exchange, with increasing demands for high-speed communication. In this context, it is crucial to guarantee a good coverage for all the users of a geographic area, but problems related to both interference \cite{mon04} and electromagnetic pollution \cite{sam17} arise in case a location is over-covered. For this reason, it is important to provide coverage to all the users, but limiting as much as possible the existence of multiple antennae covering the same area. This can clearly be modelled naturally in terms of a SCP-CS problem, making the problem we investigate important in a society which is becoming increasingly aware of sustainability \cite{wan13}.

In this paper, a mixed integer linear programming model for the SCP-CS is considered and solved by the open-source solver CP-SAT, which is part of the Google OR-Tools \cite{cpsat} optimization suite. Successful application of this solver on optimization problems with  characteristics similar to the problem under investigation, motivated our study \cite{md23}, \cite{cor}, \cite{rm25}. An experimental campaign on the  benchmark instances previously proposed in the literature is also presented and discussed.

The overall organization of the paper can be summarized as follows. The Set Covering Problem with Conflicts on Sets is formally defined in Section \ref{desc}. Section \ref{model} introduces a mixed integer linear programming model to represent the problem. In Section \ref{exp} the model proposed is  analyzed experimentally. After a description of the benchmarks previously adopted in the literature, the approach we propose is compared with the methods that previously appeared in the literature. Conclusions are finally drawn in Section \ref{conc}.

\section{Problem Description}\label{desc}

The SCP-CS can be formally defined as follows. Let $U=\{1, 2, \dots, m\}$ be a finite set, and $S=\{S_j \subseteq U | j \in N\}$ be a given collection of subsets of $U$ with $N=\{1,2,\ldots ,n\}$. For each $j \in N$ a cost $c_j \in \mathbb{Z}^+_0$, incurred once $S_j$ is selected. Moreover, for each $i,j \in N, i < j$, a cost $d_{ij} \in \mathbb{Z}^+_0$ is paid when both $S_i$ and $S_j$ are selected at the same time. We define $D \subseteq \{\{i,j\} | i,j \in N, i\ne j\}$ the set of unordered pairs $\{i,j\}$ for which $d_{ij}>0$. The objective of the SCP-CS is to select $W \subseteq S$ such that $\cup_{S_j \in W}=U$ and the total cost associated with $W$ is minimized.

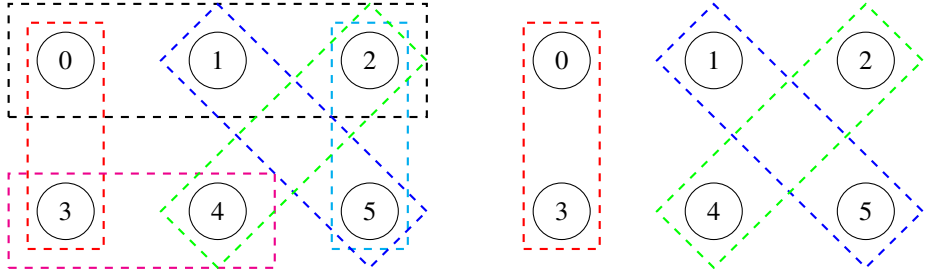
\begin{figure*}[h!]
{
\begin{center}
{
\begin{tikzpicture}[node distance={2cm}, main/.style = {draw, circle}]
			\node[main,minimum size=0.75cm] (0) {0};
			\node[white,minimum size=1.15cm] (6) [right of=0] {6};
			\node[main,minimum size=0.75cm] (1) [right of=0] {1};
			\node[white,minimum size=1.15cm] (8) [right of =1] {8};
			\node[main,minimum size=0.75cm] (2) [right of =1] {2};
			\node[main,minimum size=0.75cm] (3) [below of=0] {3};
			\node[white,minimum size=1.15cm] (9) [right of =3] {9};
			\node[main,minimum size=0.75cm] (4) [right of =3] {4};
			\node[white,main,minimum size=1.15cm] (7) [right of=4] {7};
			\node[main,minimum size=0.75cm] (5) [right of=4] {5};	
			\draw[red,thick,dashed] (-0.5,0.5) rectangle (0.5,-2.5);
			\draw[black,thick,dashed] (-0.75,0.75) rectangle (4.75,-0.75);
			\draw[magenta,thick,dashed] (-0.75,-1.5) rectangle (2.75,-2.75);
			\draw[cyan,thick,dashed] (3.5,0.5) rectangle (4.5,-2.5);
\pgfmathanglebetweenpoints{\pgfpointanchor{1}{center}}{\pgfpointanchor{5}{center}}
        \pgfmathsetmacro{\myangle}{\pgfmathresult}
        \node[draw=blue,thick,dashed, rotate fit=\myangle, fit=(6) (7)] {};
\pgfmathanglebetweenpoints{\pgfpointanchor{1}{center}}{\pgfpointanchor{5}{center}}
        \pgfmathsetmacro{\myangle}{\pgfmathresult}
        \node[draw=green,thick,dashed, rotate fit=\myangle, fit=(8) (9)] {};
\end{tikzpicture}
\hspace{1cm}
\begin{tikzpicture}[node distance={2cm}, main/.style = {draw, circle}]
			\node[main,minimum size=0.75cm] (0) {0};
			\node[white,minimum size=1.15cm] (6) [right of=0] {6};
			\node[main,minimum size=0.75cm] (1) [right of=0] {1};
			\node[white,minimum size=1.15cm] (8) [right of =1] {8};
			\node[main,minimum size=0.75cm] (2) [right of =1] {2};
			\node[main,minimum size=0.75cm] (3) [below of=0] {3};
			\node[white,minimum size=1.15cm] (9) [right of =3] {9};
			\node[main,minimum size=0.75cm] (4) [right of =3] {4};
			\node[white,main,minimum size=1.15cm] (7) [right of=4] {7};
			\node[main,minimum size=0.75cm] (5) [right of=4] {5};	
			\draw[red,thick,dashed] (-0.5,0.5) rectangle (0.5,-2.5);
\pgfmathanglebetweenpoints{\pgfpointanchor{1}{center}}{\pgfpointanchor{5}{center}}
        \pgfmathsetmacro{\myangle}{\pgfmathresult}
        \node[draw=blue,thick,dashed, rotate fit=\myangle, fit=(6) (7)] {};
\pgfmathanglebetweenpoints{\pgfpointanchor{1}{center}}{\pgfpointanchor{5}{center}}
        \pgfmathsetmacro{\myangle}{\pgfmathresult}
        \node[draw=green,thick,dashed, rotate fit=\myangle, fit=(8) (9)] {};
\end{tikzpicture}
}
	\caption{On the left an example of a simplified SCP-CS instance is presented, where the elements to be covered are numbered from 0 to 5, and the sets are the rectangular shapes, and we assume that each set has a fixed cost equal to $\alpha$. To simplify the example, no cost is shown, and we assume that two sets sharing nodes pay a penalty proportional to the overlap, and anyway much higher than the costs associated with the single sets, when they are both selected. On the right, an optimal solution with cost $3\alpha$, where three non-overlapping sets are selected to cover all the elements, and therefore no penalty is incurred. Notice that this is not always the case, and in larger problems it is common to have a penalty in the optimal cost}
	\label{figu}
\end{center}
}
\end{figure*}

A small example of an SCP-CS instance and an optimal solution for the instance are provided in Figure \ref{figu}.

\section{A Mixed Integer Linear Programming Model}\label{model}
In this section a model for the SCP-CS, obtained by extending a classic model for the Set Covering Problem \cite{vaz01} and discussed in \cite{car24}, is presented. The model is based on two sets of variables. A variable $x_{j}$ takes value 1 if set $S_j$, $ j \in N$ is selected, 0 otherwise. A variable $y_{ij}$, $\{i,j\} \in D$  takes value 1 if both sets $S_i$ and $S_j$ are selected. The resulting model is as follows.
\begin{align} 
  \min \ \ &   \sum_{j \in N} c_j x_j + \sum_{\{i,j\} \in D} d_{ij}y_{ij}& \label{1}\\ 
s.t. \ \ &	\sum_{j \in N: k \in S_j}  x_{j} \ge 1& k \in U\label{2}\\
& y_{ij} \ge x_i + x_j -1 & \{i,j\} \in D \label{3}\\
& x_{j} \in \{0,1\} & j \in N \label{4}\\
& y_{ij} \in \{0,1\} & \{i,j\} \in D \label{5}
\end{align}
The objective function (\ref{1}) minimizes the cost of the sets selected plus the penalty for conflicting selected sets.
Constraints (\ref{2}) impose that each element of set $U$ has to be covered by at least one set. 
Inequalities (\ref{3}) activate $y$ variables once the two respective conflicting sets are both selected.
Finally the domains of the variables are specified in constraints (\ref{4}) and (\ref{5}).

\section{Computational Experiments} \label{exp}
In Section \ref{ben} we describe the benchmark instances previously introduced in the literature, and used for the present study. In Section \ref{res} the approach we propose is compared with the other methods available in the literature.

\subsection{Benchmark Instances}\label{ben}
The benchmark instances for the  SCP-CS available in the literature have been introduced in \cite{car24}. They have been generated starting from the classic SCP instances made available by Beasley \cite{bea90}. Starting from a SCP instance, an instance for the SCP-CS is obtained by merging three consecutive subsets of the original instance into one (following the appearance order in the original instances), with each SCP-CS subset having a cost equal to the sum of the costs of the three original SCP subsets. This operation reduces the number of subsets to roughly one third of the original ones, but it is functional to the creation of the conflicts and the respective costs.  The value of the conflict cost $d$ is based on the overlapping of the subsets $i$ and $j$. An instance-dependent parameter $\gamma$ is calculated as follows, in order to calibrate the values of the conflicts costs:
\begin{equation}
\gamma =  \max \left \{ \left \lfloor \max_{j \in N} \left \{\frac{c_j}{|S_j|},1 \right \} \right \rceil \right\}
\end{equation}

The conflict cost for the generic subsets $i,j \in N, i<j$ is then obtained as follows:
\begin{equation}
d_{ij} = \gamma \max \left \{ |S_i \cap S_j| - \kappa, 0  \right \}
\end{equation}
where $\kappa$ is a parameter to model the number of conflicts (the higher is $\kappa$, the lower is the number of conflicts). The benchmarks proposed in \cite{car24} use $\kappa \in \{1,2\}$.

A total of  80 instances are considered. The first 50 are adapted from \cite{bea87}, 10 more instances are obtained from those introduced in \cite{gos97}, 20 from the instances originally proposed in \cite{bea92}. Considering $\kappa \in \{1,2\}$ the total number of benchmarks is 160. 
The number of elements $|U|$ ranges between 50 and 1000, the number of subsets $N$ between 167 and 3334, and the number of conflicts $|D|$ from below 2000 to more than 5.5 millions, depending on the instance and the values of the parameter $\kappa$. The detailed characteristics of the instances are summarized in Table \ref{t0}.

\begin{center}
\begin{table*}[]
\caption{Characteristics of the benchmark instances.}\label{t0}  
{\scriptsize
\centering
\begin{tabular}{lrrrrrlrrrrrlrrrr}
\toprule
Instance & $|U|$   & $|N|$  & \multicolumn{2}{c}{$|D|$}  &&Instance & $|U|$   & $|N|$  & \multicolumn{2}{c}{$|D|$}&&Instance & $|U|$   & $|N|$  & \multicolumn{2}{c}{$|D|$}   \\
 &    &   & \multicolumn{1}{c}{$\kappa = 1$ }     &  \multicolumn{1}{c}{$\kappa = 2$ }&& &    &   & \multicolumn{1}{c}{$\kappa = 1$ }     &  \multicolumn{1}{c}{$\kappa = 2$ }&& &    &   & \multicolumn{1}{c}{$\kappa = 1$ }     &  \multicolumn{1}{c}{$\kappa = 2$ }   \\
\cmidrule(lr){1-5}\cmidrule(lr){7-11}\cmidrule(lr){13-17}
scp41  & 200 & 334  & 8351   & 1908  &  & scpa3    & 300   & 1000 & 139351 & 49665  &  & scpd5   & 400  & 1334 & 883471  & 868775  \\
scp410 & 200 & 334  & 7944   & 1793  &  & scpa4    & 300   & 1000 & 139104 & 50086  &  & scpe1   & 50   & 167  & 13810   & 13707   \\
scp42  & 200 & 334  & 8306   & 1848  &  & scpa5    & 300   & 1000 & 139255 & 49635  &  & scpe2   & 50   & 167  & 13837   & 13756   \\
scp43  & 200 & 334  & 8258   & 1791  &  & scpb1    & 300   & 1000 & 485873 & 457826 &  & scpe3   & 50   & 167  & 13755   & 13666   \\
scp44  & 200 & 334  & 8596   & 1977  &  & scpb2    & 300   & 1000 & 486022 & 458275 &  & scpe4   & 50   & 167  & 13818   & 13707   \\
scp45  & 200 & 334  & 8136   & 1789  &  & scpb3    & 300   & 1000 & 485654 & 457218 &  & scpe5   & 50   & 167  & 13839   & 13790   \\
scp46  & 200 & 334  & 8826   & 2055  &  & scpb4    & 300   & 1000 & 485276 & 456577 &  & scpnre1 & 500  & 1667 & 1388611 & 1388611 \\
scp47  & 200 & 334  & 7843   & 1711  &  & scpb5    & 300   & 1000 & 485595 & 457456 &  & scpnre2 & 500  & 1667 & 1388611 & 1388611 \\
scp48  & 200 & 334  & 8526   & 1860  &  & scpc1    & 400   & 1334 & 347270 & 152214 &  & scpnre3 & 500  & 1667 & 1388611 & 1388611 \\
scp49  & 200 & 334  & 8195   & 1807  &  & scpc2    & 400   & 1334 & 345511 & 150739 &  & scpnre4 & 500  & 1667 & 1388611 & 1388611 \\
scp51  & 200 & 667  & 33962  & 8238  &  & scpc3    & 400   & 1334 & 345621 & 151469 &  & scpnre5 & 500  & 1667 & 1388611 & 1388611 \\
scp510 & 200 & 667  & 34476  & 8321  &  & scpc4    & 400   & 1334 & 346553 & 150846 &  & scpnrf1 & 500  & 1667 & 1388611 & 1388611 \\
scp52  & 200 & 667  & 34092  & 8047  &  & scpc5    & 400   & 1334 & 346137 & 151744 &  & scpnrf2 & 500  & 1667 & 1388611 & 1388611 \\
scp53  & 200 & 667  & 34153  & 8104  &  & scpclr10 & 511   & 70   & 2415   & 2415   &  & scpnrf3 & 500  & 1667 & 1388611 & 1388611 \\
scp54  & 200 & 667  & 33893  & 7913  &  & scpclr11 & 1023  & 110  & 5995   & 5995   &  & scpnrf4 & 500  & 1667 & 1388611 & 1388611 \\
scp55  & 200 & 667  & 32490  & 7376  &  & scpclr12 & 2047  & 165  & 13530  & 13530  &  & scpnrf5 & 500  & 1667 & 1388611 & 1388611 \\
scp56  & 200 & 667  & 34551  & 8355  &  & scpclr13 & 4095  & 239  & 28441  & 28441  &  & scpnrg1 & 1000 & 3334 & 4628203 & 3574676 \\
scp57  & 200 & 667  & 34832  & 8393  &  & scpcyc06 & 240   & 64   & 281    & 153    &  & scpnrg2 & 1000 & 3334 & 4626736 & 3574156 \\
scp58  & 200 & 667  & 33285  & 7784  &  & scpcyc07 & 672   & 150  & 806    & 443    &  & scpnrg3 & 1000 & 3334 & 4623974 & 3574347 \\
scp59  & 200 & 667  & 33094  & 7597  &  & scpcyc08 & 1792  & 342  & 2173   & 1230   &  & scpnrg4 & 1000 & 3334 & 4630196 & 3579650 \\
scp61  & 200 & 334  & 49416  & 40825 &  & scpcyc09 & 4608  & 768  & 5660   & 3267   &  & scpnrg5 & 1000 & 3334 & 4625519 & 3575142 \\
scp62  & 200 & 334  & 49916  & 41855 &  & scpcyc10 & 11520 & 1707 & 14370  & 8401   &  & scpnrh1 & 1000 & 3334 & 5555947 & 5555602 \\
scp63  & 200 & 334  & 49827  & 41583 &  & scpcyc11 & 28160 & 3755 & 35487  & 21072  &  & scpnrh2 & 1000 & 3334 & 5555936 & 5555589 \\
scp64  & 200 & 334  & 49404  & 40723 &  & scpd1    & 400   & 1334 & 883245 & 867920 &  & scpnrh3 & 1000 & 3334 & 5555946 & 5555566 \\
scp65  & 200 & 334  & 49830  & 41733 &  & scpd2    & 400   & 1334 & 883605 & 869019 &  & scpnrh4 & 1000 & 3334 & 5555937 & 5555574 \\
scpa1  & 300 & 1000 & 139239 & 49629 &  & scpd3    & 400   & 1334 & 883254 & 868280 &  & scpnrh5 & 1000 & 3334 & 5555925 & 5555576 \\
scpa2  & 300 & 1000 & 139174 & 49455 &  & scpd4    & 400   & 1334 & 883136 & 868196 &  &         &      &      &         &   \\
\bottomrule
\end{tabular}
}
\end{table*}
\end{center}

\subsection{Experimental Results} \label{res}
The model discussed in Section \ref{model} has been solved with the Google OR-Tools CP-SAT solver \cite{cpsat} version 9.12. The experiments have been run on a computer equipped with an Intel Core i7 12700F CPU, but the times have been normalized to the Intel
Xeon Gold 6140M 2.30 GHz CPU used for the experiments reported in \cite{car24}, in order to have the fairest possible comparison. The normalization of the computation times has been carried out according to the conversion ratio inferred from \url{http://gene.disi.unitn.it/test/cpu_list.php}. Note that a maximum computation time of 3600 seconds (on the latter processor) is allowed for all the methods. 

The methods involved in the comparison -- all originally presented in \cite{car24}, apart from the one we propose -- are:
\begin{itemize}
\item BLP: the bilinear mixed integer programming model presented in Section \ref{model} solved by Gurobi 9.5 \cite{gurobi}. The best lower and upper bound out of two runs are reported;
\item MILP: a mixed integer programming model exploiting some characteristics of the benchmarks and  solved by Gurobi 9.5. The best lower and upper bound out of two runs are reported. This model is not valid for the general SCP-CS because it exploits some special structure of the benchmark instances;
\item QBP: a bilinear quadratic programming model solved by Gurobi 9.5. The best lower and upper bound out of two runs are reported;
\item P-GRASP: a Parallel GRASP (Greedy Randomized Adaptive Search Procedure) metaheuristic algorithm. The best upper bound out of ten runs is reported;
\item CP-SAT: the mixed integer linear program presented in Section \ref{model} solved with Google OR-Tools CP-SAT solver 9.12 \cite{cpsat}. The best lower and upper bounds obtained in one run are reported.
\end{itemize}

The experimental results are summarized in Table \ref{t1} for $\kappa=1$ and in Table \ref{t2} for $\kappa=2$. For each instance the best-known lower and upper bounds from the literature are reported. Then for exact methods BLP, MILP, QBP and CP-SAT  the percentage deviations from the best known bounds are reported and the computation time in seconds required to retrieve the best known solution (\emph{Heu}) and eventually prove its optimality (\emph{Tot}). For the heuristic approach P-GRASP we report only the percentage deviation for the upper bound and the time required to retrieve this solution. The percentage deviations for lower and upper bounds are calculated respectively as $100 \cdot \frac{LB_{BK}-LB_M}{LB_{BK}}$ and $100 \cdot \frac{UB_M-UB_{BK}}{UB_{BK}}$, where the factors with subscript $BK$ are the values reported in the columns \emph{Best Known} and those with subscript $M$ are the values retrieved by the method $M$ under investigation. The last line of each table reports for each method the average values among all the instances.

The results suggest that the approach we propose positions itself as a method able to retrieve good quality lower and upper bounds, without however matching the performance achieved by the dynamic (non-compact) mixed integer linear programs in terms of lower bounds, and the quality of the upper bounds provided by specialized metaheuristics. The new approach is however able to return both lower and upper bounds of good quality, being the first method able to find such a balance, notwithstanding its simplicity. On top of this, 9 improved heuristic solutions were retrieved, They are marked in bold in the tables, and detailed in Table \ref{t2}. Finally, comparing the results of Table \ref{t1} ($\kappa=1$) with those of Table \ref{t2} ($\kappa=2$), it clearly emerges that the instances with a lower value of $\kappa$ are more difficult to solve. This is consistent with the previous literature.

\begin{center}
\begin{table}[]
\caption{Improved upper bounds.}\label{t2}  
{\scriptsize
\centering
\begin{tabular}{lrrrrrlrrrr}
\toprule
Instance & $\kappa$   &  \multicolumn{2}{c}{UB}  &&Instance & $\kappa$   & \multicolumn{2}{c}{UB}\\
 &   & Old  &New &&&    &Old &New\\
\cmidrule(lr){1-2}\cmidrule(lr){3-4}\cmidrule(lr){6-7}\cmidrule(lr){8-9}
scp54    & 1 & 1360 & 1348 &  & scpcyc08 & 2 & 659   & 645   \\
scpcyc08 & 1 & 908  & 907  &  & scpcyc09 & 2 & 1721  & 1659  \\
scpcyc10 & 1 & 5728 & 5684 &  & scpcyc10 & 2 & 4441  & 4240  \\
scpc1    & 2 & 1583 & 1518 &  & scpcyc11 & 2 & 11391 & 10476\\
\bottomrule
\end{tabular}
}
\end{table}
\end{center}

\section{Conclusions} \label{conc}
A new compact formulation based on mixed integer linear formulation for the Set Covering Problem with Conflicts on Sets was proposed and solved via an open-source solver. 

The results of the experimental campaign, run on the benchmark instances adopted in the previous literature, are reported. The results indicate that the  approach proposed is able to produce good-quality results for the instances, although the dynamic models previously appeared in the literature are able to find better lower bounds on average, and the metaheuristic approaches published before found better heuristic solutions on average. However, the new approach was able to improve the best-known upper bound (heuristic solutions) for 9 instances. 

\bibliographystyle{IEEEtran}
\bibliography{mybibfile}

\onecolumn
\begin{landscape}
{\scriptsize
\begin{longtable}{lrrrrrrrrrrrrrrrrrrrr}
\caption{Computational experiments, $\kappa = 1$.}\\  
     \toprule 
    \multicolumn{1}{c}{\multirow{4}{*}{Instances}} & \multicolumn{2}{c}{\multirow{3}{*}{Best Known}}                          & \multicolumn{4}{c}{BLP \cite{car24}}                                                                             & \multicolumn{4}{c}{MILP \cite{car24}}                                                                            & \multicolumn{4}{c}{QBP \cite{car24}}                         & \multicolumn{2}{c}{P-GRASP \cite{car24}}                       & \multicolumn{4}{c}{CP-SAT}                                                                          \\
\cmidrule(lr){4-7}\cmidrule(lr){8-11}\cmidrule(lr){12-15}\cmidrule(lr){14-17}\cmidrule(lr){18-21}
         & & & \multicolumn{2}{c}{Dev}                         & \multicolumn{2}{c}{Sec}                           & \multicolumn{2}{c}{Dev}                         & \multicolumn{2}{c}{Sec}                           & \multicolumn{2}{c}{Dev}   & \multicolumn{2}{c}{Sec}                         & \multicolumn{1}{r}{Dev} & \multicolumn{1}{r}{Sec} & \multicolumn{2}{c}{Dev}                         & \multicolumn{2}{c}{Sec}                           \\
\cmidrule(lr){2-3}\cmidrule(lr){4-5}\cmidrule(lr){6-7}\cmidrule(lr){8-9}\cmidrule(lr){10-11}\cmidrule(lr){12-13}\cmidrule(lr){14-15}\cmidrule(lr){16-16}\cmidrule(lr){17-17}\cmidrule(lr){18-19}\cmidrule(lr){20-21}

         &\multicolumn{1}{r}{LB} & \multicolumn{1}{r}{UB}   & \multicolumn{1}{r}{LB} & \multicolumn{1}{r}{UB} & \multicolumn{1}{r}{Heu} & \multicolumn{1}{r}{Tot} & \multicolumn{1}{r}{LB} & \multicolumn{1}{r}{UB} & \multicolumn{1}{r}{Heu} & \multicolumn{1}{r}{Tot} & \multicolumn{1}{r}{LB} & \multicolumn{1}{r}{UB} & \multicolumn{1}{r}{Heu} & \multicolumn{1}{r}{Tot} & \multicolumn{1}{r}{UB}  & \multicolumn{1}{r}{Heu}    & \multicolumn{1}{r}{LB} & \multicolumn{1}{r}{UB} & \multicolumn{1}{r}{Heu} & \multicolumn{1}{r}{Tot} \\
\midrule
\endfirsthead
     \caption{Computational experiments, $\kappa = 1$ (continued).}\\  
     \midrule
    \multicolumn{1}{c}{\multirow{4}{*}{Instances}} & \multicolumn{2}{c}{\multirow{3}{*}{Best Known}}                          & \multicolumn{4}{c}{BLP \cite{car24}}                                                                             & \multicolumn{4}{c}{MILP \cite{car24}}                                                                            & \multicolumn{4}{c}{QBP \cite{car24}}                         & \multicolumn{2}{c}{P-GRASP \cite{car24}}                       & \multicolumn{4}{c}{CP-SAT}                                                                          \\
\cmidrule(lr){4-7}\cmidrule(lr){8-11}\cmidrule(lr){12-15}\cmidrule(lr){14-17}\cmidrule(lr){18-21}
         & & & \multicolumn{2}{c}{Dev}                         & \multicolumn{2}{c}{Sec}                           & \multicolumn{2}{c}{Dev}                         & \multicolumn{2}{c}{Sec}                           & \multicolumn{2}{c}{Dev}   & \multicolumn{2}{c}{Sec}                         & \multicolumn{1}{r}{Dev} & \multicolumn{1}{r}{Sec} & \multicolumn{2}{c}{Dev}                         & \multicolumn{2}{c}{Sec}                           \\
\cmidrule(lr){2-3}\cmidrule(lr){4-5}\cmidrule(lr){6-7}\cmidrule(lr){8-9}\cmidrule(lr){10-11}\cmidrule(lr){12-13}\cmidrule(lr){14-15}\cmidrule(lr){16-16}\cmidrule(lr){17-17}\cmidrule(lr){18-19}\cmidrule(lr){20-21}

         &\multicolumn{1}{r}{LB} & \multicolumn{1}{r}{UB}   & \multicolumn{1}{r}{LB} & \multicolumn{1}{r}{UB} & \multicolumn{1}{r}{Heu} & \multicolumn{1}{r}{Tot} & \multicolumn{1}{r}{LB} & \multicolumn{1}{r}{UB} & \multicolumn{1}{r}{Heu} & \multicolumn{1}{r}{Tot} & \multicolumn{1}{r}{LB} & \multicolumn{1}{r}{UB} & \multicolumn{1}{r}{Heu} & \multicolumn{1}{r}{Tot} & \multicolumn{1}{r}{UB}  & \multicolumn{1}{r}{Heu}    & \multicolumn{1}{r}{LB} & \multicolumn{1}{r}{UB} & \multicolumn{1}{r}{Heu} & \multicolumn{1}{r}{Tot} \\
\midrule
\endhead
\midrule
\multicolumn{21}{r}{\textit{(Continued on next page)}} 
\endfoot
\endlastfoot
\label{t1}
scp41    & 2037 & 2037  & 0.0  & 0.0   & 241.8  & 253.7  & 0.0  & 0.0   & 199.5  & 200.3  & 0.0  & 0.0  & 109.3  & 188.8  & 2.9 & 68.8  & 23.7 & 0.0           & 170.7  & 3600.0 \\
scp410   & 2501 & 2501  & 0.0  & 0.0   & 730.5  & 903.2  & 0.0  & 0.0   & 813.9  & 911.9  & 0.0  & 0.0  & 483.6  & 511.9  & 0.0 & 103.0 & 28.3 & 0.0           & 1087.7 & 3600.0 \\
scp42    & 1977 & 1977  & 0.0  & 0.0   & 7.6    & 81.8   & 0.0  & 0.0   & 7.8    & 84.4   & 0.0  & 0.0  & 1.1    & 65.9   & 0.0 & 66.8  & 0.0  & 0.0           & 23.3   & 2581.2 \\
scp43    & 2583 & 2583  & 0.0  & 0.0   & 938.6  & 1314.3 & 0.0  & 0.0   & 1103.4 & 1428.6 & 0.0  & 0.0  & 519.5  & 2742.0 & 0.0 & 87.8  & 31.4 & 2.0           & 167.6  & 3600.0 \\
scp44    & 2543 & 2543  & 0.0  & 0.0   & 918.6  & 1033.9 & 0.0  & 0.0   & 728.1  & 824.5  & 0.0  & 0.0  & 653.7  & 662.0  & 0.0 & 166.5 & 30.9 & 0.0           & 1191.9 & 3600.0 \\
scp45    & 2247 & 2247  & 0.0  & 0.0   & 396.0  & 400.1  & 0.0  & 0.0   & 353.3  & 363.6  & 0.0  & 0.0  & 166.3  & 236.1  & 0.0 & 74.9  & 25.4 & 0.0           & 247.1  & 3600.0 \\
scp46    & 2602 & 2602  & 0.0  & 0.0   & 898.2  & 1550.4 & 0.0  & 0.0   & 901.1  & 1535.8 & 0.0  & 0.0  & 11.7   & 1081.5 & 0.0 & 91.5  & 30.5 & 2.7           & 207.8  & 3600.0 \\
scp47    & 2128 & 2128  & 0.0  & 0.0   & 747.2  & 778.8  & 0.0  & 0.0   & 693.2  & 731.5  & 0.0  & 0.0  & 1668.2 & 1860.6 & 0.0 & 2.0   & 27.8 & 0.0           & 2055.3 & 3600.0 \\
scp48    & 2647 & 2647  & 0.0  & 0.0   & 2048.3 & 2220.4 & 0.0  & 0.0   & 940.3  & 1756.2 & 0.0  & 0.0  & 1238.8 & 1496.7 & 0.6 & 84.7  & 34.5 & 0.6           & 395.0  & 3600.0 \\
scp49    & 2604 & 2604  & 0.0  & 0.0   & 268.8  & 705.3  & 0.0  & 0.0   & 759.5  & 968.8  & 13.6 & 0.0  & 1558.6 & 3600.0 & 0.0 & 16.1  & 30.8 & 0.0           & 2427.1 & 3600.0 \\
scp51    & 1305 & 1532  & 2.9  & 0.0   & 598.1  & 3600.0 & 5.1  & 5.4   & 1877.7 & 3600.0 & 0.0  & 6.9  & 3428.9 & 3600.0 & 0.0 & 199.4 & 37.2 & 0.0           & 3165.0 & 3600.0 \\
scp510   & 1338 & 1338  & 9.0  & 1.7   & 3343.5 & 3600.0 & 8.1  & 15.6  & 3203.7 & 3600.0 & 0.0  & 0.0  & 968.6  & 1291.4 & 7.2 & 95.7  & 33.3 & 0.0           & 1648.6 & 3600.0 \\
scp52    & 1446 & 1446  & 8.2  & 0.0   & 2343.8 & 3600.0 & 0.0  & 0.0   & 2927.3 & 3600.0 & 0.0  & 0.0  & 2187.3 & 2603.0 & 0.0 & 30.7  & 41.1 & 3.5           & 2805.6 & 3600.0 \\
scp53    & 1423 & 1423  & 14.4 & 5.8   & 359.8  & 3600.0 & 9.5  & 0.0   & 611.9  & 3600.0 & 0.0  & 0.0  & 2415.8 & 2864.7 & 3.9 & 152.0 & 44.3 & 0.0           & 22.0   & 3600.0 \\
scp54    & 1241 & 1360  & 8.3  & 2.3   & 93.3   & 3600.0 & 1.3  & 0.0   & 1627.8 & 3600.0 & 0.0  & 2.4  & 3321.8 & 3600.0 & 0.0 & 89.9  & 39.7 & \textbf{-0.9} & 986.8  & 3600.0 \\
scp55    & 1375 & 1375  & 2.0  & 0.0   & 3248.6 & 3600.0 & 3.6  & 0.0   & 2888.7 & 3600.0 & 0.0  & 0.0  & 2027.8 & 2085.1 & 0.5 & 143.1 & 54.7 & 5.4           & 1774.5 & 3600.0 \\
scp56    & 1217 & 1410  & 0.0  & 0.0   & 1584.0 & 3600.0 & 0.1  & 0.0   & 3256.3 & 3600.0 & 0.1  & 0.0  & 3215.3 & 3600.0 & 0.0 & 173.0 & 51.4 & 0.0           & 3526.3 & 3600.0 \\
scp57    & 1502 & 1551  & 3.5  & 0.0   & 2445.0 & 3600.0 & 5.1  & 0.0   & 2427.6 & 3600.0 & 0.0  & 0.0  & 2583.2 & 3600.0 & 1.3 & 200.7 & 52.7 & 3.9           & 2393.9 & 3600.0 \\
scp58    & 1219 & 1517  & 3.0  & 3.8   & 3339.9 & 3600.0 & 4.6  & 6.3   & 2783.1 & 3600.0 & 0.0  & 6.2  & 2976.2 & 3600.0 & 0.0 & 129.9 & 47.6 & 1.5           & 607.8  & 3600.0 \\
scp59    & 1336 & 1518  & 1.6  & 0.0   & 2110.8 & 3600.0 & 0.0  & 1.1   & 2997.2 & 3600.0 & 2.4  & 0.0  & 2637.9 & 3600.0 & 2.6 & 216.1 & 32.5 & 0.0           & 202.2  & 3600.0 \\
scp61    & 1571 & 2930  & 29.2 & 64.9  & 1791.9 & 3600.0 & 29.6 & 65.4  & 839.9  & 3600.0 & 0.0  & 0.0  & 169.1  & 3600.0 & 0.0 & 0.7   & 49.5 & 17.6          & 1912.8 & 3600.0 \\
scp62    & 1526 & 3361  & 26.5 & 43.1  & 1038.7 & 3600.0 & 26.7 & 28.2  & 1393.7 & 3600.0 & 0.0  & 2.3  & 369.6  & 3600.0 & 0.0 & 37.5  & 49.8 & 12.2          & 65.1   & 3600.0 \\
scp63    & 1634 & 4105  & 21.6 & 39.1  & 1075.2 & 3600.0 & 21.4 & 51.7  & 1041.8 & 3600.0 & 0.0  & 17.9 & 35.1   & 3600.0 & 0.0 & 40.0  & 49.0 & 23.1          & 108.3  & 3600.0 \\
scp64    & 1462 & 3211  & 25.3 & 26.5  & 2553.0 & 3600.0 & 25.3 & 52.8  & 2518.5 & 3600.0 & 0.0  & 7.3  & 1050.5 & 3600.0 & 0.0 & 30.3  & 47.3 & 20.0          & 2199.9 & 3600.0 \\
scp65    & 1364 & 3181  & 18.4 & 41.7  & 2111.7 & 3600.0 & 18.5 & 45.7  & 2609.5 & 3600.0 & 0.0  & 5.2  & 1808.3 & 3600.0 & 0.0 & 54.8  & 40.3 & 19.7          & 51.4   & 3600.0 \\
scpa1    & 1028 & 2893  & 16.4 & 28.5  & 3542.4 & 3600.0 & 16.4 & 65.4  & 3528.4 & 3600.0 & 0.0  & 8.4  & 357.6  & 3600.0 & 0.0 & 45.3  & 31.4 & 7.0           & 1114.3 & 3600.0 \\
scpa2    & 993  & 2934  & 9.8  & 37.6  & 1776.6 & 3600.0 & 9.4  & 49.8  & 3430.8 & 3600.0 & 0.0  & 12.2 & 3478.0 & 3600.0 & 0.0 & 338.5 & 27.0 & 23.4          & 1296.7 & 3600.0 \\
scpa3    & 853  & 2647  & 4.3  & 30.7  & 3234.8 & 3600.0 & 4.3  & 58.3  & 3163.0 & 3600.0 & 0.0  & 40.8 & 1998.2 & 3600.0 & 0.0 & 362.0 & 21.6 & 12.8          & 2961.4 & 3600.0 \\
scpa4    & 928  & 2701  & 10.5 & 49.0  & 3421.6 & 3600.0 & 10.5 & 57.1  & 3426.8 & 3600.0 & 0.0  & 33.3 & 887.0  & 3600.0 & 0.0 & 539.4 & 23.8 & 26.2          & 3542.4 & 3600.0 \\
scpa5    & 912  & 3069  & 4.9  & 37.4  & 2299.6 & 3600.0 & 4.9  & 62.9  & 2195.0 & 3600.0 & 0.0  & 28.1 & 115.5  & 3600.0 & 0.0 & 455.7 & 23.7 & 12.4          & 1916.0 & 3600.0 \\
scpb1    & 800  & 4320  & 0.5  & 117.1 & 10.5   & 3600.0 & 0.0  & 151.3 & 12.0   & 3600.0 & 46.1 & 15.1 & 1952.8 & 3600.0 & 0.0 & 150.9 & 61.1 & 49.0          & 2049.6 & 3600.0 \\
scpb2    & 801  & 4810  & 0.0  & 120.0 & 10.7   & 3600.0 & 0.0  & 120.0 & 12.3   & 3600.0 & 50.3 & 6.0  & 107.5  & 3600.0 & 0.0 & 198.0 & 60.0 & 29.2          & 1709.7 & 3600.0 \\
scpb3    & 767  & 4099  & 0.0  & 127.5 & 42.7   & 3600.0 & 0.0  & 127.5 & 49.8   & 3600.0 & 43.5 & 14.1 & 19.5   & 3600.0 & 0.0 & 159.7 & 58.0 & 54.3          & 3107.1 & 3600.0 \\
scpb4    & 759  & 3885  & 0.0  & 146.5 & 11.4   & 3600.0 & 0.0  & 137.6 & 13.2   & 3600.0 & 47.3 & 23.1 & 3387.8 & 3600.0 & 0.0 & 43.0  & 58.9 & 56.6          & 2158.6 & 3600.0 \\
scpb5    & 808  & 4542  & 0.0  & 140.5 & 11.3   & 3600.0 & 0.0  & 140.5 & 12.8   & 3600.0 & 49.5 & 12.7 & 1139.6 & 3600.0 & 0.0 & 156.4 & 60.8 & 53.9          & 1418.0 & 3600.0 \\
scpc1    & 806  & 7437  & 17.1 & 97.2  & 1882.3 & 3600.0 & 17.1 & 135.0 & 1729.3 & 3600.0 & 0.0  & 30.1 & 2205.4 & 3600.0 & 0.0 & 525.5 & 18.5 & 18.4          & 2905.3 & 3600.0 \\
scpc2    & 722  & 3995  & 13.2 & 87.1  & 2039.0 & 3600.0 & 13.2 & 87.1  & 2036.9 & 3600.0 & 0.0  & 12.8 & 1085.1 & 3600.0 & 0.0 & 88.9  & 17.0 & 23.1          & 2794.0 & 3600.0 \\
scpc3    & 774  & 4611  & 15.4 & 129.0 & 3008.0 & 3600.0 & 15.4 & 132.0 & 2707.8 & 3600.0 & 0.0  & 25.1 & 3539.9 & 3600.0 & 0.0 & 453.9 & 17.1 & 14.3          & 1466.0 & 3600.0 \\
scpc4    & 742  & 4457  & 17.8 & 79.4  & 3600.6 & 3600.0 & 17.8 & 79.4  & 3600.7 & 3600.0 & 0.0  & 29.9 & 246.3  & 3600.0 & 0.0 & 394.4 & 18.7 & 5.5           & 3090.6 & 3600.0 \\
scpc5    & 759  & 5239  & 18.8 & 101.5 & 1339.0 & 3600.0 & 18.8 & 101.5 & 1296.1 & 3600.0 & 0.0  & 8.3  & 1921.6 & 3600.0 & 0.0 & 355.8 & 19.8 & 9.6           & 2786.4 & 3600.0 \\
scpclr10 & 1926 & 1926  & 0.0  & 0.0   & 401.5  & 412.1  & 0.0  & 0.0   & 191.1  & 271.5  & 0.0  & 0.0  & 816.0  & 818.3  & 0.0 & 1.0   & 32.2 & 0.0           & 19.3   & 3600.0 \\
scpclr11 & 2056 & 3501  & 0.0  & 20.4  & 2065.2 & 3600.0 & 10.5 & 23.7  & 1307.1 & 3600.0 & 25.0 & 4.2  & 3116.5 & 3600.0 & 0.0 & 1.0   & 25.0 & 0.0           & 46.6   & 3600.0 \\
scpclr12 & 3103 & 6429  & 0.0  & 43.7  & 1732.0 & 3600.0 & 0.1  & 55.3  & 1838.9 & 3600.0 & 51.6 & 16.0 & 218.0  & 3600.0 & 0.0 & 21.3  & 38.0 & 10.1          & 1364.2 & 3600.0 \\
scpclr13 & 4199 & 13267 & 1.4  & 37.0  & 662.2  & 3600.0 & 0.0  & 35.1  & 211.9  & 3600.0 & 99.4 & 53.4 & 2472.9 & 3600.0 & 0.0 & 48.6  & 75.8 & 17.7          & 3155.2 & 3600.0 \\
scpcyc06 & 126  & 126   & 0.0  & 0.0   & 0.7    & 1.2    & 0.0  & 0.0   & 2.3    & 2.4    & 0.0  & 0.0  & 9.1    & 9.5    & 0.0 & 0.5   & 0.0  & 0.0           & 0.0    & 11.5   \\
scpcyc07 & 335  & 335   & 0.0  & 0.0   & 4.1    & 2273.5 & 0.0  & 0.0   & 44.8   & 1942.8 & 4.8  & 0.0  & 128.4  & 3600.0 & 0.0 & 1.2   & 20.9 & 0.0           & 21.7   & 3600.0 \\
scpcyc08 & 668  & 908   & 0.0  & 0.3   & 273.4  & 3600.0 & 0.9  & 1.2   & 56.2   & 3600.0 & 8.4  & 4.3  & 1804.5 & 3600.0 & 0.0 & 6.8   & 22.2 & \textbf{-0.1} & 108.2  & 3600.0 \\
scpcyc09 & 1369 & 2317  & 7.2  & 0.6   & 1095.0 & 3600.0 & 0.0  & 12.0  & 3600.9 & 3600.0 & 1.8  & 14.5 & 3449.3 & 3600.0 & 0.0 & 402.5 & 15.0 & 0.9           & 111.9  & 3600.0 \\
scpcyc10 & 2920 & 5728  & 0.0  & 17.9  & 3343.7 & 3600.0 & 0.3  & 32.7  & 3488.2 & 3600.0 & 9.5  & 12.1 & 3497.2 & 3600.0 & 0.0 & 456.1 & 19.6 & \textbf{-0.8} & 502.1  & 3600.0 \\
scpcyc11 & 5237 & 13727 & 1.5  & 41.2  & 0.0    & 3600.0 & 2.4  & 43.8  & 1080.6 & 3600.0 & 0.0  & 74.4 & 3324.4 & 3600.0 & 0.0 & 569.4 & 31.7 & 0.6           & 3224.8 & 3600.0 \\
scpd1    & 649  & 5340  & 0.0  & 125.8 & 22.4   & 3600.0 & 0.0  & 136.0 & 25.5   & 3600.0 & 40.5 & 20.1 & 21.1   & 3600.0 & 0.0 & 109.9 & 73.7 & 42.4          & 1039.6 & 3600.0 \\
scpd2    & 653  & 5520  & 0.0  & 111.6 & 24.1   & 3600.0 & 0.0  & 111.6 & 27.1   & 3600.0 & 43.3 & 12.2 & 3385.7 & 3600.0 & 0.0 & 263.4 & 72.9 & 26.4          & 2333.4 & 3600.0 \\
scpd3    & 651  & 5572  & 0.0  & 98.3  & 24.2   & 3600.0 & 0.0  & 98.3  & 27.2   & 3600.0 & 40.7 & 13.5 & 301.4  & 3600.0 & 0.0 & 29.0  & 80.0 & 31.8          & 686.5  & 3600.0 \\
scpd4    & 792  & 7873  & 0.0  & 71.9  & 23.1   & 3600.0 & 0.0  & 71.9  & 26.3   & 3600.0 & 53.5 & 26.2 & 975.6  & 3600.0 & 0.0 & 81.2  & 78.8 & 30.0          & 3066.5 & 3600.0 \\
scpd5    & 672  & 5718  & 0.0  & 130.4 & 26.0   & 3600.0 & 0.0  & 130.4 & 28.9   & 3600.0 & 39.6 & 15.2 & 316.4  & 3600.0 & 0.0 & 73.0  & 75.4 & 37.6          & 1306.6 & 3600.0 \\
scpe1    & 26   & 26    & 0.0  & 0.0   & 63.8   & 72.4   & 0.0  & 0.0   & 61.4   & 71.2   & 0.0  & 0.0  & 5.3    & 5.5    & 0.0 & 1.5   & 0.0  & 0.0           & 2.2    & 16.5   \\
scpe2    & 28   & 28    & 0.0  & 0.0   & 33.6   & 94.5   & 0.0  & 0.0   & 31.1   & 93.2   & 0.0  & 0.0  & 6.6    & 6.6    & 0.0 & 0.1   & 0.0  & 0.0           & 6.2    & 17.8   \\
scpe3    & 24   & 24    & 0.0  & 0.0   & 2.2    & 33.9   & 0.0  & 0.0   & 2.2    & 23.8   & 0.0  & 0.0  & 4.3    & 5.6    & 0.0 & 0.1   & 0.0  & 0.0           & 0.8    & 11.3   \\
scpe4    & 25   & 25    & 0.0  & 0.0   & 42.7   & 88.2   & 0.0  & 0.0   & 42.7   & 88.9   & 0.0  & 0.0  & 6.7    & 7.0    & 0.0 & 0.1   & 0.0  & 0.0           & 22.0   & 32.3   \\
scpe5    & 28   & 28    & 0.0  & 0.0   & 65.9   & 85.8   & 0.0  & 0.0   & 67.0   & 85.6   & 0.0  & 0.0  & 5.0    & 5.2    & 0.0 & 0.4   & 0.0  & 0.0           & 15.7   & 25.3   \\
scpnre1  & 323  & 3711  & 0.0  & 100.1 & 200.8  & 3600.0 & 0.0  & 14.1  & 224.2  & 3600.0 & 65.6 & 17.7 & 832.6  & 3600.0 & 0.0 & 333.3 & 74.0 & 43.2          & 2889.4 & 3600.0 \\
scpnre2  & 320  & 3693  & 0.0  & 145.2 & 211.1  & 3600.0 & 0.9  & 32.8  & 218.0  & 2079.3 & 69.4 & 16.9 & 1680.3 & 3600.0 & 0.0 & 263.2 & 74.7 & 39.6          & 2385.0 & 3600.0 \\
scpnre3  & 307  & 3826  & 0.0  & 114.6 & 616.8  & 3600.0 & 0.3  & 22.3  & 46.2   & 3600.0 & 63.5 & 28.4 & 1896.3 & 3600.0 & 0.0 & 137.6 & 73.6 & 35.3          & 918.6  & 3600.0 \\
scpnre4  & 305  & 3726  & 0.0  & 94.0  & 657.1  & 3600.0 & 0.0  & 13.3  & 708.9  & 3600.0 & 62.6 & 25.9 & 1780.4 & 3600.0 & 0.0 & 157.6 & 47.9 & 56.3          & 3550.0 & 3600.0 \\
scpnre5  & 329  & 3660  & 0.0  & 104.6 & 37.2   & 3600.0 & 0.0  & 40.4  & 41.7   & 2787.6 & 64.1 & 28.3 & 738.9  & 3600.0 & 0.0 & 297.4 & 52.6 & 32.5          & 2447.8 & 3600.0 \\
scpnrf1  & 157  & 1185  & 0.0  & 81.0  & 287.4  & 3600.0 & 0.0  & 54.9  & 293.7  & 3600.0 & 68.8 & 36.8 & 673.6  & 3600.0 & 0.0 & 19.0  & 41.4 & 40.2          & 1789.1 & 3600.0 \\
scpnrf2  & 165  & 1243  & 0.0  & 87.0  & 601.8  & 3600.0 & 0.6  & 51.9  & 41.1   & 3600.0 & 54.5 & 29.3 & 343.2  & 3600.0 & 0.0 & 34.2  & 45.5 & 48.1          & 389.7  & 3600.0 \\
scpnrf3  & 160  & 1184  & 0.0  & 89.3  & 35.5   & 3600.0 & 0.0  & 53.1  & 40.4   & 3600.0 & 21.3 & 32.1 & 2055.8 & 3600.0 & 0.0 & 1.8   & 41.9 & 44.3          & 1603.6 & 3600.0 \\
scpnrf4  & 162  & 1183  & 1.9  & 110.0 & 164.4  & 3600.0 & 0.0  & 43.4  & 153.7  & 3600.0 & 24.1 & 35.4 & 1568.5 & 3600.0 & 0.0 & 2.2   & 44.4 & 47.0          & 1494.3 & 3600.0 \\
scpnrf5  & 169  & 1226  & 0.0  & 49.2  & 183.3  & 3600.0 & 2.4  & 17.9  & 259.0  & 3600.0 & 66.3 & 29.7 & 2647.6 & 3600.0 & 0.0 & 78.0  & 45.6 & 40.2          & 2596.1 & 3600.0 \\
scpnrg1  & 609  & 13234 & 0.0  & 151.4 & 2860.6 & 3600.0 & 0.0  & 151.4 & 2873.1 & 3600.0 & 35.6 & 0.0  & 2712.8 & 3600.0 & 0.5 & 240.2 & 58.3 & 78.2          & 1522.1 & 3600.0 \\
scpnrg2  & 602  & 13110 & 0.0  & 165.7 & 234.4  & 3600.0 & 4.3  & 160.0 & 253.7  & 3600.0 & 38.2 & 4.2  & 241.7  & 3600.0 & 0.0 & 48.3  & 60.6 & 100.3         & 1534.0 & 3600.0 \\
scpnrg3  & 621  & 14287 & 0.0  & 97.2  & 250.5  & 3600.0 & 0.0  & 97.2  & 277.8  & 3600.0 & 34.5 & 3.1  & 156.0  & 3600.0 & 0.0 & 65.9  & 58.8 & 74.7          & 2162.6 & 3600.0 \\
scpnrg4  & 596  & 13673 & 0.0  & 190.1 & 238.5  & 3600.0 & 0.0  & 122.5 & 253.7  & 3600.0 & 33.7 & 2.8  & 220.6  & 3600.0 & 0.0 & 38.6  & 56.5 & 71.1          & 1633.1 & 3600.0 \\
scpnrg5  & 583  & 13529 & 0.0  & 99.1  & 240.3  & 3600.0 & 0.0  & 99.1  & 265.2  & 3600.0 & 33.4 & 11.4 & 2770.3 & 3600.0 & 0.0 & 90.8  & 86.3 & 78.5          & 641.3  & 3600.0 \\
scpnrh1  & 304  & 9315  & 0.7  & 107.3 & 226.7  & 3600.0 & 0.0  & 94.5  & 253.1  & 3600.0 & 60.2 & 31.5 & 2933.5 & 3600.0 & 0.0 & 129.9 & 83.2 & 72.3          & 1466.2 & 3600.0 \\
scpnrh2  & 308  & 9101  & 0.0  & 103.7 & 400.4  & 3600.0 & 1.0  & 103.7 & 786.7  & 3600.0 & 59.1 & 42.5 & 216.2  & 3600.0 & 0.0 & 14.5  & 83.4 & 101.7         & 1448.0 & 3600.0 \\
scpnrh3  & 297  & 9294  & 0.0  & 93.5  & 257.5  & 3600.0 & 0.0  & 93.5  & 274.4  & 3600.0 & 59.6 & 32.9 & 3246.0 & 3600.0 & 0.0 & 12.1  & 82.8 & 96.0          & 1462.5 & 3600.0 \\
scpnrh4  & 313  & 9644  & 0.0  & 73.7  & 892.7  & 3600.0 & 0.0  & 73.7  & 917.7  & 3600.0 & 54.0 & 24.4 & 3591.2 & 3600.0 & 0.0 & 100.9 & 83.7 & 70.1          & 1279.2 & 3600.0 \\
scpnrh5  & 300  & 9461  & 0.0  & 83.5  & 243.5  & 3600.0 & 0.0  & 65.5  & 265.6  & 3600.0 & 63.7 & 31.1 & 2860.2 & 3600.0 & 0.0 & 80.3  & 83.0 & 71.9          & 563.7  & 3600.0 \\ \midrule
Averages &      &       & 3.9  & 55.0  & 1002.6 & 2943.8 & 3.9  & 49.2  & 1091.6 & 2903.2 & 21.3 & 13.5 & 1438.5 & 2931.9 & 0.2 & 135.5 & 42.4 & 24.3          & 1432.1 & 3318.7\\
\bottomrule
\end{longtable}
}
\newpage
{\scriptsize
\begin{longtable}{lrrrrrrrrrrrrrrrrrrrr}
\caption{Computational experiments, $\kappa = 2$.}\\  
     \toprule 
    \multicolumn{1}{c}{\multirow{4}{*}{Instances}} & \multicolumn{2}{c}{\multirow{3}{*}{Best Known}}                          & \multicolumn{4}{c}{BLP \cite{car24}}                                                                             & \multicolumn{4}{c}{MILP \cite{car24}}                                                                            & \multicolumn{4}{c}{QBP \cite{car24}}                         & \multicolumn{2}{c}{P-GRASP \cite{car24}}                       & \multicolumn{4}{c}{CP-SAT}                                                                          \\
\cmidrule(lr){4-7}\cmidrule(lr){8-11}\cmidrule(lr){12-15}\cmidrule(lr){14-17}\cmidrule(lr){18-21}
         & & & \multicolumn{2}{c}{Dev}                         & \multicolumn{2}{c}{Sec}                           & \multicolumn{2}{c}{Dev}                         & \multicolumn{2}{c}{Sec}                           & \multicolumn{2}{c}{Dev}   & \multicolumn{2}{c}{Sec}                         & \multicolumn{1}{r}{Dev} & \multicolumn{1}{r}{Sec} & \multicolumn{2}{c}{Dev}                         & \multicolumn{2}{c}{Sec}                           \\
\cmidrule(lr){2-3}\cmidrule(lr){4-5}\cmidrule(lr){6-7}\cmidrule(lr){8-9}\cmidrule(lr){10-11}\cmidrule(lr){12-13}\cmidrule(lr){14-15}\cmidrule(lr){16-16}\cmidrule(lr){17-17}\cmidrule(lr){18-19}\cmidrule(lr){20-21}

         &\multicolumn{1}{r}{LB} & \multicolumn{1}{r}{UB}   & \multicolumn{1}{r}{LB} & \multicolumn{1}{r}{UB} & \multicolumn{1}{r}{Heu} & \multicolumn{1}{r}{Tot} & \multicolumn{1}{r}{LB} & \multicolumn{1}{r}{UB} & \multicolumn{1}{r}{Heu} & \multicolumn{1}{r}{Tot} & \multicolumn{1}{r}{LB} & \multicolumn{1}{r}{UB} & \multicolumn{1}{r}{Heu} & \multicolumn{1}{r}{Tot} & \multicolumn{1}{r}{UB}  & \multicolumn{1}{r}{Heu}    & \multicolumn{1}{r}{LB} & \multicolumn{1}{r}{UB} & \multicolumn{1}{r}{Heu} & \multicolumn{1}{r}{Tot} \\
\midrule
\endfirsthead
     \caption{Computational experiments, $\kappa = 2$ (continued).}\\  
     \midrule
    \multicolumn{1}{c}{\multirow{4}{*}{Instances}} & \multicolumn{2}{c}{\multirow{3}{*}{Best Known}}                          & \multicolumn{4}{c}{BLP \cite{car24}}                                                                             & \multicolumn{4}{c}{MILP \cite{car24}}                                                                            & \multicolumn{4}{c}{QBP \cite{car24}}                         & \multicolumn{2}{c}{P-GRASP \cite{car24}}                       & \multicolumn{4}{c}{CP-SAT}                                                                          \\
\cmidrule(lr){4-7}\cmidrule(lr){8-11}\cmidrule(lr){12-15}\cmidrule(lr){14-17}\cmidrule(lr){18-21}
         & & & \multicolumn{2}{c}{Dev}                         & \multicolumn{2}{c}{Sec}                           & \multicolumn{2}{c}{Dev}                         & \multicolumn{2}{c}{Sec}                           & \multicolumn{2}{c}{Dev}   & \multicolumn{2}{c}{Sec}                         & \multicolumn{1}{r}{Dev} & \multicolumn{1}{r}{Sec} & \multicolumn{2}{c}{Dev}                         & \multicolumn{2}{c}{Sec}                           \\
\cmidrule(lr){2-3}\cmidrule(lr){4-5}\cmidrule(lr){6-7}\cmidrule(lr){8-9}\cmidrule(lr){10-11}\cmidrule(lr){12-13}\cmidrule(lr){14-15}\cmidrule(lr){16-16}\cmidrule(lr){17-17}\cmidrule(lr){18-19}\cmidrule(lr){20-21}

         &\multicolumn{1}{r}{LB} & \multicolumn{1}{r}{UB}   & \multicolumn{1}{r}{LB} & \multicolumn{1}{r}{UB} & \multicolumn{1}{r}{Heu} & \multicolumn{1}{r}{Tot} & \multicolumn{1}{r}{LB} & \multicolumn{1}{r}{UB} & \multicolumn{1}{r}{Heu} & \multicolumn{1}{r}{Tot} & \multicolumn{1}{r}{LB} & \multicolumn{1}{r}{UB} & \multicolumn{1}{r}{Heu} & \multicolumn{1}{r}{Tot} & \multicolumn{1}{r}{UB}  & \multicolumn{1}{r}{Heu}    & \multicolumn{1}{r}{LB} & \multicolumn{1}{r}{UB} & \multicolumn{1}{r}{Heu} & \multicolumn{1}{r}{Tot} \\
\midrule
\endhead
\midrule
\multicolumn{21}{r}{\textit{(Continued on next page)}} 
\endfoot
\endlastfoot
\label{t2}
scp41    & 1108 & 1108  & 0.0  & 0.0   & 1.5    & 1.5    & 0.0  & 0.0      & 1.4    & 1.5    & 0.0  & 0.0  & 1.9    & 1.9    & 0.0 & 3.3   & 0.0  & 0.0           & 0.6    & 1.5    \\
scp410   & 1404 & 1404  & 0.0  & 0.0   & 2.4    & 2.8    & 0.0  & 0.0      & 2.1    & 2.3    & 0.0  & 0.0  & 1.8    & 2.1    & 0.0 & 11.2  & 0.0  & 0.0           & 0.3    & 2.9    \\
scp42    & 1209 & 1209  & 0.0  & 0.0   & 0.5    & 22.5   & 0.0  & 0.0      & 0.5    & 1.4    & 0.0  & 0.0  & 0.5    & 2.2    & 0.0 & 4.5   & 0.0  & 0.0           & 2.4    & 3.3    \\
scp43    & 1113 & 1113  & 0.0  & 0.0   & 0.3    & 0.9    & 0.0  & 0.0      & 0.3    & 0.9    & 0.0  & 0.0  & 0.2    & 0.9    & 0.0 & 1.1   & 0.0  & 0.0           & 0.5    & 1.4    \\
scp44    & 1192 & 1192  & 0.0  & 0.0   & 2.5    & 2.6    & 0.0  & 0.0      & 2.5    & 2.6    & 0.0  & 0.0  & 2.0    & 2.2    & 0.0 & 3.0   & 0.0  & 0.0           & 1.1    & 1.9    \\
scp45    & 1279 & 1279  & 0.0  & 0.0   & 1.8    & 2.0    & 0.0  & 0.0      & 1.7    & 1.8    & 0.0  & 0.0  & 1.5    & 1.7    & 0.0 & 51.0  & 0.0  & 0.0           & 2.4    & 3.1    \\
scp46    & 1302 & 1302  & 0.0  & 0.0   & 1.6    & 3.2    & 0.0  & 0.0      & 1.6    & 3.2    & 0.0  & 0.0  & 1.6    & 3.1    & 0.0 & 60.6  & 0.0  & 0.0           & 0.3    & 1.8    \\
scp47    & 1116 & 1116  & 0.0  & 0.0   & 2.2    & 2.4    & 0.0  & 0.0      & 2.1    & 2.3    & 0.0  & 0.0  & 2.7    & 2.8    & 0.0 & 5.1   & 0.0  & 0.0           & 0.5    & 1.3    \\
scp48    & 1149 & 1149  & 0.0  & 0.0   & 0.3    & 1.7    & 0.0  & 0.0      & 0.4    & 1.5    & 0.0  & 0.0  & 1.3    & 1.4    & 0.0 & 8.5   & 0.0  & 0.0           & 0.2    & 1.1    \\
scp49    & 1398 & 1398  & 0.0  & 0.0   & 4.5    & 4.8    & 0.0  & 0.0      & 4.5    & 4.8    & 0.0  & 0.0  & 1.3    & 2.3    & 0.0 & 53.3  & 0.0  & 0.0           & 1.2    & 2.5    \\
scp51    & 618  & 618   & 0.0  & 0.0   & 2.7    & 4.9    & 0.0  & 0.0      & 2.8    & 5.4    & 0.0  & 0.0  & 1.7    & 1.9    & 0.0 & 13.5  & 0.0  & 0.0           & 3.7    & 5.8    \\
scp510   & 642  & 642   & 0.0  & 0.0   & 0.8    & 1.9    & 0.0  & 0.0      & 0.9    & 1.9    & 0.0  & 0.0  & 0.7    & 0.9    & 0.0 & 23.3  & 0.0  & 0.0           & 0.9    & 2.3    \\
scp52    & 602  & 602   & 0.0  & 0.0   & 0.2    & 1.2    & 0.0  & 0.0      & 0.2    & 1.2    & 0.0  & 0.0  & 0.4    & 0.9    & 0.0 & 5.2   & 0.0  & 0.0           & 0.5    & 1.6    \\
scp53    & 627  & 627   & 0.0  & 0.0   & 4.7    & 5.1    & 0.0  & 0.0      & 4.9    & 5.4    & 0.0  & 0.0  & 1.6    & 1.8    & 0.0 & 14.6  & 0.0  & 0.0           & 2.4    & 6.8    \\
scp54    & 546  & 546   & 0.0  & 0.0   & 2.1    & 2.1    & 0.0  & 0.0      & 2.1    & 2.1    & 0.0  & 0.0  & 0.8    & 0.9    & 0.0 & 75.8  & 0.0  & 0.0           & 1.2    & 2.3    \\
scp55    & 528  & 528   & 0.0  & 0.0   & 0.6    & 0.8    & 0.0  & 0.0      & 0.7    & 0.8    & 0.0  & 0.0  & 0.2    & 0.4    & 0.0 & 2.9   & 0.0  & 0.0           & 0.5    & 1.0    \\
scp56    & 511  & 511   & 0.0  & 0.0   & 0.5    & 0.6    & 0.0  & 0.0      & 0.4    & 0.5    & 0.0  & 0.0  & 0.3    & 0.3    & 0.0 & 6.5   & 0.0  & 0.0           & 0.4    & 0.7    \\
scp57    & 764  & 764   & 0.0  & 0.0   & 7.2    & 7.6    & 0.0  & 0.0      & 7.4    & 7.7    & 0.0  & 0.0  & 2.0    & 2.1    & 0.0 & 77.4  & 0.0  & 0.0           & 1.8    & 6.0    \\
scp58    & 650  & 650   & 0.0  & 0.0   & 1.1    & 3.5    & 0.0  & 0.0      & 1.1    & 3.3    & 0.0  & 0.0  & 0.6    & 1.4    & 0.0 & 140.9 & 0.0  & 0.0           & 1.2    & 6.3    \\
scp59    & 660  & 660   & 0.0  & 0.0   & 0.9    & 3.2    & 0.0  & 0.0      & 1.0    & 3.3    & 0.0  & 0.0  & 1.0    & 1.0    & 0.0 & 8.6   & 0.0  & 0.0           & 1.2    & 2.6    \\
scp61    & 1375 & 1733  & 24.4 & 25.2  & 2805.3 & 3600.0 & 26.3 & 31.9     & 3383.3 & 3600.0 & 0.0  & 4.3  & 696.0  & 3600.0 & 0.0 & 31.2  & 50.9 & 4.3           & 3252.3 & 3600.0 \\
scp62    & 1260 & 2037  & 31.7 & 23.6  & 160.1  & 3600.0 & 31.7 & 34.5     & 177.4  & 3600.0 & 0.0  & 5.4  & 208.3  & 3600.0 & 0.0 & 28.3  & 47.9 & 0.0           & 627.0  & 3600.0 \\
scp63    & 1430 & 2355  & 18.0 & 24.2  & 1417.2 & 3600.0 & 18.0 & 31.8     & 1406.6 & 3600.0 & 0.0  & 6.1  & 952.4  & 3600.0 & 0.0 & 25.5  & 48.7 & 0.0           & 2592.8 & 3600.0 \\
scp64    & 1299 & 1899  & 21.1 & 12.1  & 1477.9 & 3600.0 & 36.9 & 53.0     & 2886.9 & 3600.0 & 0.0  & 0.0  & 2314.3 & 3600.0 & 0.0 & 33.8  & 49.7 & 8.9           & 336.2  & 3600.0 \\
scp65    & 1385 & 1960  & 32.2 & 18.3  & 2489.0 & 3600.0 & 31.6 & 32.1     & 3493.8 & 3600.0 & 0.0  & 1.7  & 629.8  & 3600.0 & 0.0 & 43.2  & 49.9 & 1.7           & 64.2   & 3600.0 \\
scpa1    & 906  & 963   & 0.0  & 0.0   & 1019.5 & 3600.0 & 0.9  & 0.0      & 1100.1 & 3600.0 & 0.2  & 1.3  & 2917.2 & 3600.0 & 1.6 & 343.2 & 26.3 & 1.6           & 3438.3 & 3600.0 \\
scpa2    & 890  & 1048  & 0.0  & 0.0   & 3417.0 & 3600.0 & 0.0  & 0.0      & 2205.8 & 3600.0 & 1.2  & 7.5  & 2540.6 & 3600.0 & 0.5 & 141.4 & 21.7 & 0.0           & 1278.3 & 3600.0 \\
scpa3    & 875  & 907   & 1.6  & 0.3   & 8.7    & 3600.0 & 1.6  & 0.3      & 9.7    & 3600.0 & 0.0  & 0.0  & 1997.6 & 3600.0 & 0.0 & 204.4 & 26.9 & 0.0           & 2007.7 & 3600.0 \\
scpa4    & 946  & 946   & 7.4  & 0.0   & 311.6  & 3600.0 & 5.1  & 0.0      & 206.5  & 3600.0 & 0.0  & 0.0  & 2775.2 & 3238.9 & 0.0 & 257.1 & 25.5 & 0.0           & 2243.7 & 3600.0 \\
scpa5    & 894  & 894   & 0.0  & 0.0   & 1868.4 & 2010.3 & 0.0  & 0.0      & 40.4   & 697.6  & 0.0  & 0.0  & 255.9  & 1576.3 & 0.0 & 43.0  & 22.5 & 3.1           & 1683.1 & 3600.0 \\
scpb1    & 596  & 2863  & 0.0  & 111.1 & 10.5   & 3600.0 & 0.0  & 133.9    & 15.4   & 3600.0 & 25.2 & 3.9  & 2966.0 & 3600.0 & 0.0 & 290.6 & 43.3 & 13.8          & 2960.8 & 3600.0 \\
scpb2    & 618  & 3044  & 0.0  & 100.7 & 1548.1 & 3600.0 & 0.0  & 100.7    & 1566.6 & 3600.0 & 21.2 & 20.3 & 2203.9 & 3600.0 & 0.0 & 138.1 & 46.3 & 34.5          & 3363.9 & 3600.0 \\
scpb3    & 571  & 2643  & 0.0  & 121.9 & 11.0   & 3600.0 & 0.0  & 81.7     & 12.3   & 3600.0 & 19.6 & 17.1 & 2858.0 & 3600.0 & 0.0 & 246.6 & 41.2 & 34.4          & 923.6  & 3600.0 \\
scpb4    & 574  & 2470  & 0.0  & 108.3 & 10.0   & 3600.0 & 0.0  & 108.3    & 11.4   & 3600.0 & 27.0 & 26.7 & 95.9   & 3600.0 & 0.0 & 88.5  & 41.1 & 24.0          & 487.9  & 3600.0 \\
scpb5    & 611  & 2944  & 0.0  & 109.9 & 1642.1 & 3600.0 & 0.0  & 109.9    & 1660.3 & 3600.0 & 28.8 & 8.2  & 685.0  & 3600.0 & 0.0 & 264.2 & 42.4 & 17.2          & 272.5  & 3600.0 \\
scpc1    & 733  & 1583  & 0.0  & 1.9   & 581.4  & 3600.0 & 1.2  & 0.1      & 638.4  & 3600.0 & 2.2  & 11.6 & 138.8  & 3600.0 & 0.0 & 573.2 & 11.1 & \textbf{-4.1} & 2311.7 & 3600.0 \\
scpc2    & 706  & 1249  & 8.8  & 11.4  & 2257.4 & 3600.0 & 8.8  & 27.8     & 2244.2 & 3600.0 & 0.0  & 17.7 & 944.7  & 3600.0 & 0.0 & 423.4 & 18.0 & 8.9           & 1273.2 & 3600.0 \\
scpc3    & 693  & 1343  & 11.1 & 17.9  & 1340.3 & 3600.0 & 11.1 & 15.6     & 1322.7 & 3600.0 & 0.0  & 30.0 & 199.3  & 3600.0 & 0.0 & 484.0 & 11.1 & 6.3           & 910.6  & 3600.0 \\
scpc4    & 699  & 1385  & 9.4  & 15.2  & 1786.3 & 3600.0 & 9.4  & 18.1     & 1784.1 & 3600.0 & 0.0  & 11.4 & 3147.2 & 3600.0 & 0.0 & 488.6 & 13.3 & 8.5           & 670.9  & 3600.0 \\
scpc5    & 720  & 1295  & 9.6  & 16.9  & 1606.9 & 3600.0 & 9.6  & 5.8      & 1583.3 & 3600.0 & 0.0  & 18.4 & 1715.5 & 3600.0 & 0.0 & 281.6 & 17.6 & 0.6           & 2746.4 & 3600.0 \\
scpclr10 & 1871 & 1871  & 0.0  & 0.0   & 253.1  & 331.9  & 0.0  & 0.0 & 184.5  & 258.8  & 0.0  & 0.0  & 476.2  & 809.0  & 0.0 & 0.2   & 0.0  & 0.0           & 10.1   & 1102.4 \\
scpclr11 & 1797 & 3446  & 0.0  & 17.1  & 2659.2 & 3600.0 & 0.1  & 0.0      & 2951.7 & 3600.0 & 8.0  & 0.0  & 1502.6 & 3600.0 & 0.0 & 0.8   & 6.8  & 0.0           & 13.6   & 3600.0 \\
scpclr12 & 3076 & 6374  & 0.0  & 21.2  & 2657.1 & 3600.0 & 0.1  & 35.1     & 1932.9 & 3600.0 & 55.8 & 0.0  & 580.4  & 3600.0 & 0.0 & 21.7  & 22.2 & 0.0           & 278.1  & 3600.0 \\
scpclr13 & 4184 & 13212 & 1.3  & 106.8 & 0.0    & 3600.0 & 0.0  & 61.2     & 1072.3 & 3600.0 & 96.4 & 59.8 & 270.7  & 3600.0 & 0.0 & 3.6   & 70.2 & 11.6          & 232.0  & 3600.0 \\
scpcyc06 & 99   & 99    & 0.0  & 0.0   & 0.5    & 0.8    & 0.0  & 0.0      & 0.4    & 0.7    & 0.0  & 0.0  & 7.9    & 8.0    & 0.0 & 0.6   & 0.0  & 0.0           & 0.1    & 2.9    \\
scpcyc07 & 250  & 250   & 0.0  & 0.0   & 147.4  & 1688.1 & 0.0  & 0.0      & 2.8    & 1246.7 & 2.4  & 0.0  & 436.8  & 3600.0 & 0.0 & 33.0  & 13.2 & 0.0           & 1.1    & 3600.0 \\
scpcyc08 & 507  & 659   & 0.0  & 0.2   & 208.8  & 3600.0 & 0.0  & 0.0      & 2112.6 & 3600.0 & 3.2  & 0.6  & 475.7  & 3600.0 & 2.7 & 366.0 & 11.8 & \textbf{-2.1} & 2216.6 & 3600.0 \\
scpcyc09 & 1080 & 1721  & 0.0  & -2.5  & 2454.3 & 3600.0 & 2.4  & 0.0      & 1924.3 & 3600.0 & 5.1  & 2.4  & 3230.4 & 3600.0 & 1.7 & 448.9 & 11.2 & \textbf{-3.6} & 115.7  & 3600.0 \\
scpcyc10 & 2216 & 4441  & 0.2  & 0.5   & 3600.1 & 3600.0 & 0.0  & 0.0      & 3595.4 & 3600.0 & 3.4  & 3.4  & 2989.1 & 3600.0 & 1.9 & 569.0 & 5.6  & \textbf{-4.5} & 1707.7 & 3600.0 \\
scpcyc11 & 4765 & 11391 & 0.4  & 11.4  & 2740.5 & 3600.0 & 0.4  & 20.4     & 2771.7 & 3600.0 & 0.0  & 26.3 & 3285.8 & 3600.0 & 0.0 & 288.6 & 14.5 & \textbf{-8.0} & 2924.3 & 3600.0 \\
scpd1    & 605  & 4104  & 0.0  & 116.4 & 23.9   & 3600.0 & 0.0  & 116.4    & 27.0   & 3600.0 & 44.3 & 37.3 & 1877.8 & 3600.0 & 0.0 & 317.5 & 56.2 & 29.0          & 2790.0 & 3600.0 \\
scpd2    & 579  & 4213  & 0.0  & 119.7 & 24.7   & 3600.0 & 0.0  & 119.7    & 27.5   & 3600.0 & 36.1 & 6.5  & 477.6  & 3600.0 & 0.0 & 276.2 & 53.9 & 32.2          & 2741.5 & 3600.0 \\
scpd3    & 596  & 4468  & 0.0  & 90.9  & 23.2   & 3600.0 & 0.0  & 105.4    & 26.0   & 3600.0 & 39.1 & 9.4  & 316.1  & 3600.0 & 0.0 & 89.4  & 54.2 & 28.2          & 2613.2 & 3600.0 \\
scpd4    & 670  & 6209  & 0.0  & 121.3 & 23.2   & 3600.0 & 0.0  & 121.3    & 26.1   & 3600.0 & 50.4 & 14.4 & 985.6  & 3600.0 & 0.0 & 509.0 & 58.2 & 20.0          & 657.3  & 3600.0 \\
scpd5    & 595  & 4287  & 0.0  & 96.0  & 24.1   & 3600.0 & 0.0  & 96.0     & 27.2   & 3600.0 & 42.9 & 10.6 & 2171.4 & 3600.0 & 0.0 & 412.7 & 52.8 & 17.4          & 3530.7 & 3600.0 \\
scpe1    & 23   & 23    & 0.0  & 0.0   & 30.0   & 69.6   & 0.0  & 0.0      & 20.8   & 50.6   & 0.0  & 0.0  & 4.6    & 6.5    & 0.0 & 0.9   & 0.0  & 0.0           & 13.4   & 21.0   \\
scpe2    & 25   & 23    & 0.0  & 0.0   & 17.6   & 63.1   & 0.0  & 0.0      & 15.8   & 71.4   & 0.0  & 0.0  & 4.1    & 10.3   & 0.0 & 0.5   & 0.0  & 8.7           & 2.3    & 29.3   \\
scpe3    & 23   & 22    & 0.0  & 0.0   & 8.2    & 42.4   & 0.0  & 0.0      & 7.1    & 54.1   & 0.0  & 0.0  & 5.6    & 9.2    & 0.0 & 0.1   & 0.0  & 4.5           & 0.3    & 20.3   \\
scpe4    & 22   & 22    & 0.0  & 0.0   & 8.6    & 30.2   & 0.0  & 0.0      & 15.6   & 45.3   & 0.0  & 0.0  & 4.1    & 9.9    & 0.0 & 0.1   & 0.0  & 0.0           & 2.6    & 17.1   \\
scpe5    & 25   & 25    & 0.0  & 0.0   & 20.5   & 73.8   & 0.0  & 0.0      & 20.6   & 73.8   & 0.0  & 0.0  & 6.2    & 6.7    & 0.0 & 0.1   & 0.0  & 0.0           & 1.5    & 22.4   \\
scpnre1  & 317  & 3648  & 0.0  & 93.2  & 930.3  & 3600.0 & 0.0  & 30.8     & 973.9  & 3600.0 & 60.3 & 17.8 & 1033.3 & 3600.0 & 0.0 & 5.2   & 68.1 & 29.3          & 2071.6 & 3600.0 \\
scpnre2  & 312  & 3449  & 0.0  & 119.0 & 206.6  & 3600.0 & 0.0  & 51.1     & 231.4  & 1546.8 & 61.5 & 29.7 & 1963.8 & 3600.0 & 0.0 & 7.8   & 73.7 & 38.9          & 586.1  & 3600.0 \\
scpnre3  & 309  & 3549  & 0.0  & 85.4  & 3451.7 & 3600.0 & 0.0  & 58.5     & 3600.3 & 2072.9 & 66.3 & 29.5 & 2734.3 & 3600.0 & 0.0 & 465.8 & 57.9 & 30.8          & 1985.2 & 3600.0 \\
scpnre4  & 307  & 3558  & 0.0  & 77.8  & 706.1  & 3600.0 & 0.0  & 36.5     & 618.6  & 3600.0 & 63.5 & 15.1 & 364.6  & 3600.0 & 0.0 & 145.6 & 74.6 & 72.8          & 3368.3 & 3600.0 \\
scpnre5  & 311  & 3579  & 0.0  & 104.6 & 37.2   & 3600.0 & 0.0  & 26.6     & 41.8   & 3600.0 & 65.3 & 26.7 & 383.4  & 3600.0 & 0.0 & 123.3 & 74.9 & 41.2          & 526.5  & 3600.0 \\
scpnrf1  & 159  & 1175  & 0.0  & 138.1 & 33.4   & 3600.0 & 0.0  & 53.8     & 38.5   & 1360.8 & 74.8 & 31.6 & 2203.4 & 3600.0 & 0.0 & 78.5  & 42.8 & 60.4          & 836.8  & 3600.0 \\
scpnrf2  & 162  & 1227  & 0.0  & 80.4  & 248.0  & 3600.0 & 0.0  & 60.4     & 258.0  & 3600.0 & 61.1 & 30.6 & 1777.8 & 3600.0 & 0.0 & 484.0 & 40.1 & 16.5          & 395.7  & 3600.0 \\
scpnrf3  & 165  & 1174  & 0.0  & 77.3  & 867.5  & 3600.0 & 0.0  & 61.2     & 744.5  & 2459.3 & 29.1 & 32.5 & 1058.5 & 3600.0 & 0.0 & 41.3  & 30.9 & 58.2          & 3500.7 & 3600.0 \\
scpnrf4  & 161  & 1173  & 0.0  & 88.1  & 455.2  & 3600.0 & 3.1  & 55.9     & 168.4  & 1921.1 & 32.9 & 30.7 & 2735.9 & 3600.0 & 0.0 & 357.6 & 19.3 & 49.6          & 2263.8 & 3600.0 \\
scpnrf5  & 161  & 1216  & 0.0  & 49.2  & 589.0  & 3600.0 & 0.0  & 40.5     & 137.0  & 1765.8 & 54.7 & 11.3 & 1857.6 & 3600.0 & 0.0 & 204.6 & 43.5 & 48.6          & 3572.5 & 3600.0 \\
scpnrg1  & 486  & 6110  & 0.0  & 157.6 & 443.5  & 3600.0 & 0.0  & 255.7    & 458.4  & 3600.0 & 28.4 & 0.0  & 2453.9 & 3600.0 & 0.7 & 392.8 & 45.7 & 89.2          & 3093.2 & 3600.0 \\
scpnrg2  & 459  & 6272  & 0.0  & 147.5 & 227.5  & 3600.0 & 0.0  & 147.5    & 238.8  & 3600.0 & 29.2 & 13.4 & 53.8   & 3600.0 & 0.0 & 575.9 & 39.4 & 82.0          & 3525.0 & 3600.0 \\
scpnrg3  & 484  & 6947  & 0.0  & 153.9 & 204.3  & 3600.0 & 0.0  & 153.9    & 217.8  & 3600.0 & 27.9 & 12.0 & 67.0   & 3600.0 & 0.0 & 587.2 & 45.5 & 76.2          & 1122.6 & 3600.0 \\
scpnrg4  & 482  & 6589  & 0.0  & 167.2 & 246.7  & 3600.0 & 0.0  & 167.2    & 248.6  & 3600.0 & 22.4 & 1.0  & 2525.5 & 3600.0 & 0.0 & 348.4 & 44.8 & 62.4          & 2499.5 & 3600.0 \\
scpnrg5  & 476  & 6493  & 0.0  & 143.7 & 1038.2 & 3600.0 & 0.0  & 143.7    & 1075.6 & 3600.0 & 28.6 & 31.2 & 2470.9 & 3600.0 & 0.0 & 344.5 & 41.0 & 70.7          & 3533.9 & 3600.0 \\
scpnrh1  & 304  & 8574  & 0.0  & 134.7 & 248.7  & 3600.0 & 0.7  & 105.7    & 287.3  & 3600.0 & 59.9 & 32.6 & 2244.8 & 3600.0 & 0.0 & 160.6 & 76.3 & 62.9          & 2038.7 & 3600.0 \\
scpnrh2  & 308  & 8880  & 0.0  & 119.9 & 386.7  & 3600.0 & 0.6  & 89.9     & 274.1  & 3600.0 & 55.2 & 22.6 & 3047.1 & 3600.0 & 0.0 & 9.1   & 83.4 & 81.7          & 1937.7 & 3600.0 \\
scpnrh3  & 294  & 8833  & 0.0  & 93.4  & 228.4  & 3600.0 & 0.0  & 93.4     & 250.5  & 3600.0 & 60.9 & 24.0 & 3249.8 & 3600.0 & 0.0 & 140.4 & 83.7 & 76.1          & 1339.8 & 3600.0 \\
scpnrh4  & 313  & 8688  & 0.0  & 93.5  & 253.6  & 3600.0 & 1.9  & 93.5     & 274.0  & 3600.0 & 62.3 & 22.0 & 3025.8 & 3600.0 & 0.0 & 117.2 & 78.3 & 71.5          & 1857.0 & 3600.0 \\
scpnrh5  & 291  & 9070  & 0.0  & 81.5  & 253.0  & 3600.0 & 0.0  & 65.0     & 269.8  & 3600.0 & 59.8 & 24.7 & 114.9  & 3600.0 & 0.0 & 103.8 & 77.7 & 73.7          & 1608.4 & 3600.0 \\ \midrule
Averages &      &       & 2.2  & 45.3  & 647.3  & 2349.8 & 2.5  & 40.6   & 662.0  & 2196.0 & 19.0 & 10.0 & 1059.3 & 2366.4 & 0.1 & 162.9 & 27.8 & 18.6          & 1187.3 & 2400.9\\
\bottomrule
\end{longtable}
}
\end{landscape}
\twocolumn

\end{document}